\documentclass[11pt,a4paper]{article}

\usepackage{hyperref}

\usepackage{amsmath}
\usepackage{amssymb}
\usepackage{amsthm}
\usepackage{amscd}
\usepackage{amsfonts}

\usepackage[osf,sc]{mathpazo}
\usepackage{euscript}

\DeclareMathOperator{\Hom}{Hom}
\DeclareMathOperator{\End}{End}

\DeclareMathOperator{\Res}{Res}
\renewcommand{\det}{\mathrm{det}}
\DeclareMathOperator{\Tr}{Tr}

\theoremstyle{plain}
\newtheorem{lemma}{Lemma}
\newtheorem{proposition}{Proposition}
\newtheorem{corollary}{Corollary}
\newtheorem{theorem}{Theorem}

\theoremstyle{remark}

\theoremstyle{definition}
\newtheorem{definition}{Definition}
\newtheorem{remark}{Remark}
\newtheorem{example}{Example}

\def\Z{\mathbf{Z}} 
\def\R{\mathbf{R}} 
\def\Q{\mathbf{Q}} 
\def\F{\mathbf{F}} 
\def\A{\mathbf{A}} 
\def\Gm{\mathbf{G}_m} 
\def\G{\mathbf{G}} 
\def\m{\mathfrak{m}} 
\def\p{\mathfrak{p}} 

\def\O{{\cal O}} 

\renewcommand{\mod}{\operatorname{mod}}
\renewcommand{\phi}{\varphi}

\DeclareMathOperator{\GL}{GL}

\newcommand{\onto}{\mbox{\mathsurround=0pt \;$\longrightarrow \hspace{-0.7em} \to$\;}}

\newcommand{\into}{ \mbox{\mathsurround=0pt \;\raisebox{0.63ex}{\small $\subset$} \hspace{-1.07em} $\longrightarrow$\;}}

\newcommand{\covol}{{\rm covol}}
\newcommand{\Lie}{{\rm Lie}}
\newcommand{\Ext}{{\rm Ext}}

\begin{document}

\author{Lenny Taelman}
\title{Special $L$-values of Drinfeld modules}

\maketitle

\begin{abstract}

We state and prove a formula for a certain value of the Goss 
$L$-function of a Drinfeld module.  
This gives characteristic-$p$-valued function field analogues of the class number formula
and of the Birch and Swinnerton-Dyer conjecture.
The formula and its proof are presented in an entirely self-contained
fashion.

{\it Mathematics Subject Classification (2000)}: 
11R58,  
11G09.  
\end{abstract}

\thispagestyle{empty}

{\small
\tableofcontents
}

\section{Introduction and statement of the theorem}\label{intro}

Let $ k $ be a finite field of $ q $ elements. 
For a finite $ k[ t ] $-module $ M $ we define $|M| \in k[T] $ to be 
the characteristic polynomial of the endomorphism $ t $ of the $ k $-vector space $ M $,
so
\[
	|M| = \det_{k[T]}\!\left( T - t\, \big| \,M \right).
\]
In other words, if $M\cong \oplus_i k[t]/f_i(t)$ with the $f_i$ monic then
$|M| = \prod_i f_i(T)$.
We will treat the invariant
$ |M| $ of the finite $ k[ t ] $-module $ M $ as an analogue of the cardinality
$ \# A $ of a finite abelian group $ A $.

Let $ R $ be the integral closure of $ k[ t ] $
in a finite extension $ K $ of the field $ k(t) $ of rational functions. Consider 
in the Laurent series field $ k((T^{-1})) $ the infinite sum 
\[
	\sum_{I} \frac{1}{|R/I|}
\]
where $ I $ ranges over all the non-zero ideals of $ R $.
Since there are only finitely many ideals of given index, the sum converges to an element of
$ 1 + T^{-1} k[[ T^{-1} ]] $ which we denote by $ \zeta( R, 1 ) $. By unique factorization into prime ideals, we also have
\[
	\zeta( R, 1 ) = \prod_{\m} \left( 1 - \frac{1}{|R/\m|} \right)^{-1}
\]
where $ \m $ ranges over all the maximal ideals of $ R $. We stress that
$\zeta(R,1)$ depends not only on the ring structure of $ R $ but also on its $k[t]$-algebra structure.

A particular case of our main result will be a formula for $ \zeta( R, 1 ) $ analogous to the class number formula for the residue at $ s=1 $ of the Dedekind zeta function of a number field.
Whereas the class number formula is essentially a statement about the multiplicative group
$ \Gm $, our formula for $ \zeta( R, 1 ) $ will 
essentially be a statement about the \emph{Carlitz module}:
 
\begin{definition}
The \emph{Carlitz module} is the functor
\[
	C\colon \{\text{ $k[t]$-algebras } \} 
	\to \{\text{ $ k[t] $-modules } \}
\]
that associates with a $ k[t] $-algebra $ B $ the $ k[t] $-module
$ C( B ) $ whose underlying $ k $-vector space is $ B $ and whose
$ k[t] $-module structure is given by the homomorphism of $ k $-algebras
\[
	\phi_C\colon k[t] \to \End_k( B )\colon t \mapsto t + \tau,
\]
where $ \tau $ denotes the Frobenius endomorphism $b \mapsto b^q$.
\end{definition}

This definition might appear rather ad hoc, yet the functor $ C $ is in many ways analogous 
to the functor
\[
	\Gm\colon \{\text{ $\Z$-algebras } \} 
	\to \{\text{ $\Z$-modules } \}.
\]
For example, in analogy with $ \# \Gm( F ) = \# F - 1 $ for a finite field $ F $, we have

\begin{proposition}\label{minusone}
If $ \ell $ is a $ k[t] $-algebra that is a finite field then $ |C( \ell )| = | \ell | - 1 $.
\end{proposition}

In particular, this proposition allows us to re-express the infinite product defining
$ \zeta( R, 1 ) $ as follows:
\[
	\zeta( R, 1 ) = \prod_{ \m  } \frac{|R/\m|}{|C(R/\m)|} 
	\in 1 + T^{-1}k[[T^{-1}]].
\]

\begin{proof}[Proof of Proposition \ref{minusone}]
We need to show that for every finite extension $ \ell/k $ and for every element $ t \in \ell $
we have
\[
	\det_{k[T]}\!\left( T - t - \tau \,\big|\, \ell \right) = 
	\det_{k[T]}\!\left( T - t \,\big|\, \ell \right) - 1.
\]
We will do so by extending scalars from
$ k $ to an algebraic closure $ \bar{k} $.
Let $ S $ be the set of $k$-algebra embeddings of $ \ell $ into $ \bar{k} $. We have an isomorphism of $ \bar{k} $-algebras
\[
	 \ell \otimes_k\bar{k}   \to \bar{k}^S\colon
	x \otimes y \mapsto ( \sigma(x)y )_{\sigma \in S}.
\]
Under this isomorphism the action of $ t \otimes 1 $ on $ \ell \otimes \bar{k} $ 
corresponds to the diagonal action
of $ (\sigma(t))_{\sigma\in S} $ on $ \bar{k}^S $, and that of $ \tau\otimes 1 $ to
a cyclic permutation of $ S $. From this description it is easy to verify that
\[
	\det_{\bar{k}[T]}\!\left(
		T - t\otimes 1 - \tau\otimes 1 \,\big|\, \ell\otimes_k \bar{k}
	\right) = 
	\det_{\bar{k}[T]}\!\left(
		T - t\otimes 1 \,\big|\, \ell\otimes_k \bar{k} \right) - 1,
\]
as desired.
\end{proof}

We will also consider objects more general than the Carlitz module:

\begin{definition}
Let $ n $ be a non-negative integer.
A \emph{Drinfeld module of rank $ n $} over $ R $ is a functor
\[
	E\colon \{\text{ $R$-algebras } \} 
	\to \{\text{ $ k[t] $-modules } \}
\]
that associates with a $ k[t] $-algebra $ B $ the $ k[t] $-module
$ E( B ) $ whose underlying $ k $-vector space is $ B $ and whose
$ k[t] $-module structure is given by 
\[
	\phi_E\colon k[t] \to \End_k(B)\colon 
	t \mapsto t + r_1\tau + \cdots + r_n\tau^n,
\]
for fixed $ r_i \in R $, depending only on $ E $, and with $ r_n \neq 0 $.
\end{definition}

\begin{remark}
Equivalently a Drinfeld module over $ R $ is a $k[t]$-module scheme $ E $ over $ R $
whose underlying $ k $-vector space scheme is $ \G_{a,R} $ and where the induced action
of $ k[t] $ on $ \Lie( E ) = R $ (the tangent space at zero) is the tautological one.
\end{remark}

\begin{remark}
The above definition is more general than the standard one in that it allows primes of bad reduction ($r_n$ is not required to be a unit), and that it allows the trivial case $n=0$.
\end{remark}

For a Drinfeld module $ E $ over $ R $ we define $ L( E/R ) $ as follows:
\[
	L( E/R ) := \prod_{\m} \frac{|R/\m|}{|E(R/\m)|}
	\in 1 + T^{-1}k[[T^{-1}]].
\]
It is not hard to show that this infinite product indeed converges, but in any case it
will follow from the proof of the main result below. 

\begin{example}
The Carlitz module $C$ is a Drinfeld module of rank $ 1 $ and  $ \zeta( R, 1 ) = L( C/R ) $.
\end{example}

Let $ K $ be the field of fractions of $ R $ and denote the tensor product
$ R \otimes_{k[t]} k((t^{-1})) $ by $ K_\infty $. Note that $ K_\infty $ is isomorphic to a product of finite extensions of the field $ k((t^{-1})) $ and that $ R $ is discrete
and co-compact in $ K_\infty $.

\begin{proposition} \label{exp} 
For every Drinfeld module $ E $
there exists a unique power series
\[
	\exp_E X  = X + e_1 X^q + e_2 X^{q^2} + e_3 X^{q^3} + \cdots \in K_\infty[[ X ]]
\]
such that 
\begin{equation}\label{functionalequation}
	\exp_E( tX ) = \phi_E(t) (\exp_E X ).
\end{equation}
The power series $ \exp_E $ has infinite radius of convergence and defines a map
\[
	\exp_E\colon K_\infty \to E( K_\infty )
\]
which is $ k[t] $-linear, continuous and open.
\end{proposition}

\begin{proof} (See also {\cite[\S 3]{Drinfeld74E}}.) The proof is straightforward: the second
condition gives a recursion relation for the coefficients $e_i$ which determines 
them uniquely, and a closer look at the recursion relation reveals that the $ e_i $ tend to zero fast enough for $ \exp_E $ to have infinite radius of convergence. Clearly $ \exp_E $
defines a continuous $k$-linear function on $ K_\infty $ and the functional equation
(\ref{functionalequation}) guarantees that it is $ k[t] $-linear when seen as a map from
$ K_\infty $ to $ E(K_\infty) $. The non-archimedean implicit function theorem 
(\emph{see} for example \cite[2.2]{Igusa00}) implies that $ \exp_E $ is open. 
\end{proof}

\begin{proposition}\label{unittheorem} 
$\exp_E^{-1} E(R) \subset K_\infty $ is discrete and co-compact.
\end{proposition}

\begin{remark}
In particular $ \exp_E^{-1} E(R) $ is a finitely generated $k[t]$-module. The module
$ E( R ) $, however, is only finitely generated in case 
$ E $ is the Drinfeld module of rank $ 0 $ (\emph{see} \cite{Poonen97}).
\end{remark}

\begin{proof}[Proof of Proposition \ref{unittheorem}] (See also \cite{Taelman10}.) This follows from the openness of $ \exp_E $
and the fact that $ E(R) $ is discrete and co-compact in $ E(K_\infty) $. 
\end{proof}

Let $ V $ be a finite-dimensional $ k((t^{-1})) $-vector space. A \emph{lattice}
in $ V $ is by definition a discrete and co-compact sub-$k[t]$-module. The set
$ \mathcal{L}(V) $ of lattices in $ V $ carries a natural topology.

\begin{proposition}\label{defreg}
There is a unique function
\[
	\left[ - : - \right] \colon
	\mathcal{L}(V) \times \mathcal{L}(V)  \to k((T^{-1}))
\]
that is  continuous in both arguments and with the property that 
if $ \Lambda_1 $ and $ \Lambda_2 $ are sub-lattices of a lattice
$ \Lambda $ in $ V $ then 
\[
	\left[ \Lambda_1 : \Lambda_2 \right] =
	\frac{ |\Lambda/\Lambda_2| }{ |\Lambda/\Lambda_1| }.
\]
\end{proposition}

\begin{proof}
Given a lattice $ \Lambda_1 $ the set of $ \Lambda_2 $ so that $ \Lambda_1 $ and
$ \Lambda_2 $ are contained in some common lattice $ \Lambda $ is dense
in $ \mathcal{L}(V) $, so if the
function exists it is necessarily unique.
Let $ \i $ denote the obvious isomorphism
\[
	k((t^{-1})) \to k((T^{-1}))\colon \sum a_i t^i \mapsto \sum a_iT^i.
\] 
Given lattices $ \Lambda_1 $ and $ \Lambda_2 $ there exists a $ \sigma \in \GL(V) $
so that $ \sigma(\Lambda_1) = \Lambda_2 $. If we define $ [\Lambda_1 : \Lambda_2 ] $ to
be the unique monic representative of $ k^\times \i(\det(\sigma)) $ (which does not
depend on the choice of $\sigma$), then it is easy to verify that this defines
a function $ [ - : - ] $ with the required properties. 
\end{proof}

We define $ H( E/R ) $ to be the $ k[ t ] $-module
\[
	H( E/R ) := \frac{ E( K_\infty ) }{ E( R) + \exp_E  K_\infty  }. 
\]
\begin{proposition} $ H( E/R ) $ is finite.
\end{proposition}

\begin{proof} (See also \cite{Taelman10}.) The quotient $ H( E/R ) $ is
compact because $E(R)$ is co-compact in $ E(K_\infty) $, and discrete
because $\exp_E K_\infty  $ is open in $ E(K_\infty ) $.
\end{proof}

In this paper we shall prove the following result.

\begin{theorem}\label{mainthm}
For every Drinfeld module $ E $ over $ R $ we have
\[
 L( E/R ) = \left[\,  R : \exp^{-1}_E E(R) \,\right] \cdot \big| H( E/R ) \big|.
\]
\end{theorem}

\begin{example} If $n=0$ then $\exp_E $ is the identity map and both left- and right-hand-side of the formula equal $1$.
\end{example}

\begin{example} The smallest nontrivial example is for $E=C$, with $R=k[t]$. One can
verify directly that $H(C/k[t])=0$ and $\exp^{-1} C(k[t])$ is of rank $1$, generated by 
$\exp^{-1}(1)$. One thus recovers the formula
\[
	\zeta( k[t], 1 ) = \exp^{-1}(1)
\]
that was already known to Carlitz \cite[p. 160]{Carlitz35} \cite[p.199]{Anderson96}.
\end{example}

\begin{remark}
The quantity $ [  R : \exp^{-1}_E E(R) ] \in k((T^{-1})) $ is typically transcendental over $ k(T) $, \emph{see} \cite{Papanikolas08}, \cite{Chang10}, \cite{Chang10b}.
\end{remark}

\begin{remark} 
Although we will not need this, we briefly indicate how $ L(E/R) $ is the value
at some integer $ s $ of a \emph{Goss $ L $-function} (\emph{see} \cite{Goss92} for
the theory of such $ L $-functions).

If the leading coefficient $ r_n $ 
of $ \phi_E(t) $ is a unit in $ R $ (in other words: if $ E $ has everywhere good reduction),
then $ L( E/R ) $ is the value at $s=0$ of the Goss $ L $-function $ L( E^\vee_K, s ) $ defined through the dual of the Tate module of $ E_K $. 
Using Gardeyn's Euler factors at places of bad reduction \cite{Gardeyn02}, one can
define $ L( E^\vee_K, s ) $ even when $ r_n $ is not a unit. Our $ L( E/R ) $ however
depends on the chosen integral model, and will in general differ from
$ L( E^\vee_K, 0 ) $ by a factor in $ k(T)^\times $.
\end{remark}

\begin{remark} The statement of Theorem \ref{mainthm} is surprisingly similar to the
statement of the class number formula. Let $ F $ be a number field and $\O_F$ be the ring of integers in $ F $. Consider
the exponential map
\[
	\exp\colon  (\O_F/\Z) \otimes_\Z \R \to (\O_F\otimes_\Z \R)^\times/\R^\times_{>0}.
\]
Both $ \O_F/\Z $ and $ \exp^{-1} \O_F^\times  $ are discrete and co-compact subgroups
of $ V := (\O_F/\Z) \otimes_\Z \R $. For the latter this follows from the Dirichlet unit theorem. The class number formula, usually printed as
\[
	\lim_{s\to 1} (s-1)\zeta_F(s) =
	\frac{ 2^{r_1} (2\pi)^{r_2} R_F h_F }{ w_F \sqrt{ |\Delta_F| } },
\]
implies that 
\[
	\lim_{s\to 1} (s-1)\zeta_F(s) = \lambda \cdot
	\frac{ \covol(  \exp^{-1}( \O_F^\times ) ) }
		{ \covol( \O_F/\Z ) } ,
\]
for some $ \lambda \in \Q^\times $, if both co-volumes are computed
with respect to the same Haar measure on $ V $
({\it see} for example \cite[{\sc i.\S 8}]{Tate84}).
\end{remark}

\begin{remark}
If $ E $ is a Drinfeld module of rank $2$ then the Euler factors in the product
defining the left-hand side
in Theorem \ref{mainthm} share  many properties with the Euler factors $ L_v( A, 1 ) $ 
for an elliptic curve $ A $ over a number field. In this case the
theorem gives a Drinfeld module analogue of the Birch and Swinnerton-Dyer
conjecture. In fact, it is heuristics for the growth of the infinite product
\[
	\prod_{p} \frac{p}{\# A(\F_p) }
\]
that led Birch and Swinnerton-Dyer to the statement of their conjecture, \emph{see}
\cite{SwinnertonDyer67}.
\end{remark}

\begin{remark} \label{tamagawa}
Let $ \A_K $ be the ring of ad\`{e}les of $ K $. With a reasonable
theory of $ k((T^{-1})) $-valued ``measures'' on topological $ k[t] $-modules
(extending our notion of ``cardinality'' $|M|$ of finite $ k[t] $-modules $M$)
one should be able to express Theorem \ref{mainthm} as a kind of Tamagawa number
statement, namely as stating that the ``volume''  of the compact $k[t]$-module
$ E( \A_K )/E( K ) $  equals the ``volume'' of the $k[t]$-module $ \A_K/K $. (See also Remark \ref{volumes}.)
\end{remark}

\begin{remark}
A generalization of Theorem \ref{mainthm} from Drinfeld modules to
Anderson's abelian $t$-modules \cite{Anderson86},
should give, amongst others, a formula for the zeta values
\[
	\zeta( R, n ) := \sum_{I} \frac{1}{|R/I|^n}
\]
for $ n>1 $, generalizing results of Anderson and Thakur \cite{Anderson90}
for $ R = k[t] $.
\end{remark}

\begin{remark}
Going even further, it should be possible to formulate and prove an equivariant 
version of Theorem \ref{mainthm}, for abelian $t$-modules equipped with an action
of a group (or an associative algebra).  This should lead to
Stark ``units'' in extensions of function fields. 
\end{remark}

\begin{remark}
We have only considered ``$ \infty $-adic'' special values, that is, values in
$ k((T^{-1})) $. In analogy with classical $ p $-adic special values,
one could ask for $ v$-adic special values associated with Drinfeld modules, with $ v $ a finite prime of $ k(T) $. Vincent Lafforgue has expressed those in terms of extensions of
shtukas and the equal characteristic analogue of Fontaine theory  \cite{Lafforgue09}.
In the special case of the Carlitz module it is shown in \cite{Taelman11} that
$ \exp_E^{-1}E(R) $ and $ H(E/R) $ can be computed as an $\Ext^1$ respectively $\Ext^2$
in the category of shtukas considered by Lafforgue. 
\end{remark}

\bigskip

Finally, we end this introduction with a short overview of the proof of Theorem \ref{mainthm}.

Consider the power series
\[
	\Theta := \frac{ 1-\phi_E(t)T^{-1} }{ 1 - tT^{-1} } - 1 = 
	\sum_{n=1}^\infty (t-\phi_E(t))t^{n-1}T^{-n} 
\]
in $T^{-1}$ and $\tau$. For all maximal ideals $\m$ of $R$ we have
\[
	\frac{|R/\m|}{|E(R/\m)|} = 
	\frac{
		\det_{k[T]}\!\left( T - t \,\big|\, R/\m \right)
	}{
		\det_{k[T]}\!\left( T - \phi_E(t) \,\big|\, R/\m \right)
	}
	=
	\det_{k[[T^{-1}]]}\!\left( 1 + \Theta \,\big|\, R/\m \right)^{-1}.
\]
Using the trace formula of \S \ref{sectraceformula} (essentially a variation on the Woods Hole fixed-point theorem) we obtain
\begin{eqnarray*}
	L( E/R ) &=& 
	\prod_\m \det_{k[[T^{-1}]]}\!\left( 1 + \Theta \,\big|\, R/\m \right)^{-1} \\
	&=&
	\det_{k[[T^{-1}]]}\!\left( 1 + \Theta \,\big|\, \frac{ K_\infty }{ R } \right),
\end{eqnarray*}
where the last determinant is an infinite-dimensional determinant which ``converges'' in 
the sense of \S \ref{secdeterminants} (one should think of $\Theta$ as a trace-class operator
and of the determinant as a Fredholm determinant). So we have replaced the infinite product of $L$-factors by a single determinant, at the expense of having to deal with $k$-vector spaces of infinite dimension. 

Now enters the exponential map. For simplicity, assume that $H(E/R)$ vanishes. Then the exponential 
power series defines a $k$-linear isomorphism
\[
	\exp_E\colon \frac{K_\infty}{\exp^{-1} E(R)} \to \frac{K_\infty}{R}
\]
and by the functional equation (\ref{functionalequation})
we have
\[
	1 + \Theta = \frac{ 1 - \exp_E t\exp_E^{-1}T^{-1} }{ 1 - tT^{-1} }  
\]
as endomorphisms of $\frac{K_\infty[[T^{-1}]]}{R[[T^{-1}]]}$.
If $K_\infty/R$ and $E(K_\infty)/E(R)$ were not just compact but actually \emph{finite} we would conclude
\[
	\det_{k[[T^{-1}]]}\!\left( 1 + \Theta \,\big|\, \frac{ K_\infty }{ R } \right)
	= \frac{ | K_\infty/\exp_E^{-1} E(R) | }{ | K_\infty/R | },
\]
but of course they are not finite and the right-hand-side does not make any sense.
However we do show in \S \ref{secratio} that in some sense the determinant can be interpreted
as the ratio of the ``volumes'' of $ K_\infty/\exp_E^{-1} E(R) $ and
$  K_\infty/R $. More precisely, we show 
\[
	\det_{k[[T^{-1}]]}\!\left( 1 + \Theta \,\big|\, \frac{ K_\infty }{ R } \right)
	= [ R : \exp_E^{-1} E(R) ].
\]
Since we assumed $|H(E/R)|=1$, this then yields the claimed identity 
$L(E/R) = |H(E/R)| \cdot [ R : \exp_E^{-1} E(R) ]$.

\bigskip

\emph{Notation and conventions}

$k$ is a finite field of $ q $ elements, fixed throughout the text.

If $ M $ is a finite $ k[t] $-module then $ |M| $ is by definition the characteristic
polynomial of the action of $ t $ on $ M $, in the variable $ T $.
Namely,
\[
	|M| := \det_{k[T]}\!\left( T - t\, \big| \,M \right) \in k[T].
\]
The cardinality of a set $ S $ will be denoted by $ \# S $. (For example,
$ \# M = q^{\deg |M| } $.) 

If $ \alpha $ and $ \beta $ are endomorphisms of some abelian group, and if $ \beta $ is invertible then
\[
	\frac{ \alpha }{\beta}
\]
will mean $ \alpha \beta^{-1} $. 

If $ F $ is a field with a non-archimedean absolute value
$ \| \cdot \|\colon F \to \mathbf{R}_{\geq 0} $
then a (non-archimedean) \emph{norm} $ \| \cdot \| $ on an $ F $-vector space $ V $ is a map
\[
	\| \cdot \| \colon V \to \mathbf{R}_{\geq 0 }
\]
satisfying
\begin{enumerate}
\item $ \| v \| = 0 $ if and only if $ v=0 $;
\item $ \| \lambda v \| = \| \lambda\|\cdot \| v \| $ for all $ v\in V $ and $ \lambda \in F $;
\item $ \| v + w \| \leq \max( \|v\|, \|w\| ) $.
\end{enumerate}
A norm defines a metric and a topology on $ V $. We equip the finite field $ k $ with the trivial absolute value ($ \| \lambda \| = 1 $ for all $ \lambda \neq 0 $) and the
field $ k((t^{-1})) $ with an absolute value determining
the natural topology on $ k((t^{-1})) $. In particular we have $ \| t^{-1} \| < 1 $. 

\section{Nuclear operators and determinants} \label{secdeterminants}

We will make use of determinants of certain class of endomorphisms of infinite-dimensional
vector spaces and will use some terminology inspired by functional analysis for this.
It must be stressed however, that 
this is mostly a bookkeeping tool, all ``functional analysis'' here is 
merely a reflection of linear algebra on finite-dimensional
quotients of the spaces considered. 

The reader that is familiar with Anderson's ``trace calculus'' \cite[\S 2]{Anderson00}
(or Tate's traces of ``finite potent'' endomorphisms \cite[\S 1]{Tate68}) will easily
recognize that the contents of this section is mostly an adaptation of their work. We should
warn the same reader, however, that our notion of nucleus is dual to that of Anderson. (See also Remark \ref{andersontraceformula} in \S \ref{sectraceformula}.)

\bigskip

Let $ V =  ( V, \|\cdot \| ) $ be a normed $ k $-vector space and
let $ \phi $ be a continuous endomorphism of $ V $.

\begin{definition}
We say that $ \phi $ is \emph{locally contracting} if there is an
open subspace $ U \subset V $ and a real number  $ 0<c<1 $ so that
\[
	\| \phi(v) \| \leq  c\| v \|
\]
for all $ v $ in $ U $. Any such open subspace $ U $
which moreover satisfies 
$ \phi(U) \subset U $ is called a \emph{nucleus} for $ \phi $. An
open subspace $ U \subset V $ with $ \phi(U) \subset U $ and so that
$ \phi^n(U) $ is a nucleus for some $ n > 0 $ is called a \emph{pre-nucleus}.
\end{definition}

\begin{example} 
Let $ F $ be a local field containing $ k $ and let $ \| \cdot \| $ be an absolute value
for $ F $. Let $ \O_F $ denote the ring of integers. Then $ \| \cdot \| $ is a norm on
the $ k $-vector space $ \O_F $ and the continuous endomorphism
\[
	\tau\colon \O_F \to \O_F\colon f \mapsto f^q
\]
is locally contracting. The maximal ideal $ \m_F \subset \O_F $ is a nucleus for $ \tau $.
\end{example}

\begin{example} With the notation of the previous example. 
Let $ R \subset F $ be a discrete and co-compact sub-$k$-algebra. Then $ \| \cdot \| $
induces a norm on the $ k $-vector space $ F/R $. The continuous endomorphism
\[
	\tau\colon \frac{F}{R} \to \frac{F}{R} \colon f \mapsto f^q
\]
is locally contracting and the image of $ \m_F $ in $ F/R $ is a nucleus.
\end{example}

\begin{proposition}\label{commonnucleus}
Any finite collection of locally contracting endomorphisms of 
$ V $ has a common nucleus.
\end{proposition}

\begin{proof}
It suffices to observe that if $ \phi $ is a locally contracting
endomorphism then there is an $ \epsilon>0$ 
so that the open ball of radius $ \epsilon $ around the origin is a nucleus for $ \phi $.
\end{proof}

\begin{proposition}\label{sumprod}
If $ \phi $ and $ \psi $ are locally contracting, then so are the sum $ \phi + \psi $ and
composition $ \phi \psi $. 
\end{proposition}

\begin{proof} Clear.
\end{proof}

Let $ V $ be a normed $ k $-vector space. For every positive integer $ N $ we
denote by $ V[[Z]]/Z^N $ the $ k[[Z]/Z^N $-module$ V \otimes_k k[[Z]]/Z^N $
and by $ V[[Z]] $ the $k[[Z]]$-module $ V[[Z]] := \varprojlim_n V[[Z]]/Z^N $ 
equipped with the limit topology.  

Any continuous $ k[[Z]] $-linear endomorphism 
\[
	 \Phi\colon V[[Z]] \to V[[Z]]
\] 
is of the form
\[
	\Phi = \sum_{n=0}^\infty \phi_nZ^n
\]
where the $ \phi_i $ are continuous endomorphisms of $ V $. Similarly, any continuous
$ k[[Z]]/Z^n $-linear endomorphism of $ V[[Z]]/Z^N $ is of the form 
\[
	\Phi = \sum_{n=0}^{N-1} \phi_n Z^n.
\]

From now one, we assume that
$ V $ is \emph{compact}. In particular every open subspace $ U \subset V $ is of finite
co-dimension.

\begin{definition}
We say that the continuous $ k[[Z]] $-linear endomorphism
$ \Phi $ of $ V[[Z]] $ (resp. of $V[[Z]]/Z^N$)
is \emph{nuclear} if for all $n $ (resp. for all $n<N$) the 
endomorphism $\phi_n$ of $ V $ is locally contracting. 
\end{definition}

\begin{proposition}\label{existencedet}
Assume that $ \Phi\colon V[[Z]]/Z^N \to V[[Z]]/Z^N $ is nuclear.
Let $ U_1 $ and $ U_2 $ be common nuclei for the $ \phi_n $ with $ n<N $. Then 
\[
	\det_{k[[Z]]/Z^N}\!\left( 1 + \Phi\, \big|\, (V/U_i)[[Z]]/Z^N \right) 
	\in k[[Z]]/Z^N
\]
is independent of $ i\in\{1,2\} $.
\end{proposition}

\begin{proof}
It suffices to show that if $ U' \subset U $ are common nuclei for the $ \phi_n $ with
$ n < N $ then 
\[
	\det_{k[[Z]]/Z^N}\!\left( 1 + \Phi \,\big|\, (U/U')\otimes_k k[[Z]]/Z^N \right)
	= 1.
\]
Intersecting $ U $ with open balls of varying radius we can find
\[
	U = U_0 \supset U_1 \supset \cdots \supset U_m = U'
\]
with $ \phi_n( U_i ) \subset U_{i+1} $ for all $ n<N $ and $i<m $. Clearly for every $ i $
we have
\[
	\det_{k[[Z]]/Z^N}\!\left( 1 + \Phi \,\big|\, (U_i/U_{i+1}) \otimes_k k[[Z]]/Z^N \right)
	= 1,
\]
and the proposition follows.
\end{proof}

\begin{definition} If $ \Phi $ is a nuclear endomorphism of $ V[[Z]]/Z^N $ then
we denote the determinant in Proposition \ref{existencedet} by
\[
	\det_{k[[Z]]/Z^N}\! \left( 1 + \Phi \,\big|\, V \right).
\]
If $ \Phi $ is a nuclear endomorphism of $ V[[Z]] $ then we denote by 
\[
	\det_{k[[Z]]}\! \left( 1 + \Phi \,\big|\, V \right) \in k[[Z]]
\]
the unique power series that reduces to 
\[
	\det_{k[[Z]]/Z^N}\! \left( 1 + \Phi \,\big|\, V \right)
\]
modulo $ Z^N $ for every $ N $.
\end{definition}

\begin{remark}
It would be more correct to denote the above determinants by
\[
	\det_{k[[Z]]/Z^N}\! \left( 1 + \Phi \,\big|\, V[[Z]]/Z^N \right)
\]
and
\[
	\det_{k[[Z]]}\! \left( 1 + \Phi \,\big|\, V[[Z]] \right),
\]
but we will generally drop the ``$[[Z]]$'' and ``$[[Z]]/Z^N$'' from $ V $ 
in order not to overload the notation.	
\end{remark}

\begin{example}
If $ V $ is finite-dimensional over $k$ then any continuous $k[[Z]]$-linear endomorphism $ \Phi $ of $ V[[Z]]$ is
nuclear with nucleus $ 0 \subset V $ and
\[
	\det_{k[[Z]]}\! \left( 1 + \Phi \,\big|\, V \right)
\]
coincides with the determinant in the usual sense.
\end{example}

\begin{example}
If $ \phi\colon V \to V $ is locally contracting 
then we call 
\[
	\det_{k[[Z]]}\!\left( 1 - \phi Z \,\big|\, V \right) \in k[[Z]]
\]
the \emph{characteristic power series} of $ \phi $ (which is in fact a polynomial).
We have
\[
	\det_{k[[Z]]}\!\left( 1 - \phi Z \,\big|\, V \right)
	=
	\det_{k[[Z]]}\!\left( 1 - \phi Z \,\big|\, V/U \right)
\]
for any pre-nucleus $ U $.
\end{example}

\begin{proposition}[Multiplicativity in short exact sequences]
\label{seqmultiplicativity}
Let $ \Phi $ be a nuclear endomorphism of $ V[[Z]] $. Let $ W \subset V $ be
a closed subspace so that
\[
	 \Phi( W[[Z]] ) \subset W[[Z]]. 
\]
Then $ \Phi $ is nuclear on $ W[[Z]] $ and $ (V/W)[[Z]] $ and
\[
	\det_{k[[Z]]}\!\left( 1 + \Phi \,\big|\, V \right) =
	\det_{k[[Z]]}\!\left( 1 + \Phi \,\big|\, W \right)
	\det_{k[[Z]]}\!\left( 1 + \Phi \,\big|\, V/W \right).
\]
\end{proposition}

\begin{proof} Clear from the multiplicativity of finite determinants.
\end{proof}

\begin{proposition}[Multiplicativity for composition]
\label{endmultiplicativity}
Let $ \Phi $ and $ \Psi $ be nuclear endomorphisms of $ V[[Z]] $. 
Then $ ( 1 +\Phi)(1+\Psi) - 1 $ is nuclear and
\[
	\det_{k[[Z]]}\!\left( (1 + \Phi)(1 + \Psi) \,\big|\, V \right) =
	\det_{k[[Z]]}\!\left( 1 + \Phi \,\big|\, V \right) 
	\det_{k[[Z]]}\!\left( 1 + \Psi \,\big|\, V \right).
\] 
\end{proposition}

\begin{proof} Clear from the multiplicativity of finite determinants.
\end{proof}

\begin{theorem}\label{commutators}
If $ \phi $, $ \phi \alpha $, and $ \alpha\phi $ are locally contracting then
\[
	\det_{k[[Z]]}\!\left(  1 - \phi\alpha Z \,\big|\, V \right) = 
	\det_{k[[Z]]}\!\left(  1 - \alpha\phi Z \,\big|\, V \right).
\]
\end{theorem}

\begin{proof}
(See also Anderson \cite[Prop. 8]{Anderson00} and Tate \cite[(T${}_5$)]{Tate68}.)
Let $ U $ be a common nucleus for $ \alpha\phi $, $ \phi\alpha $ and $ \phi $. Consider the iterated inverse images
\[
	U_n := (\alpha\phi)^{-n}( U ) = \{ u \in V \,\big|\, (\alpha\phi)^n(u) \in U \}. 
\]
Clearly the $ U_n $ are pre-nuclei for $ \alpha\phi $. We claim that for every non-negative integer
$ n $ the subspace $ \alpha^{-1}(U_n) $ is a pre-nucleus for $ \phi\alpha $.
Indeed, since $\alpha\phi( U_n) \subset U_n $ we have the implication
\[
	\alpha(v) \in U_n \implies \alpha\phi\alpha(v) \in U_n
\]
hence
\[
	v \in \alpha^{-1}( U_n ) \implies \phi\alpha(v) \in \alpha^{-1}( U_n ).
\]
Moreover, we have
\[
	(\phi\alpha)^{n+1}(\alpha^{-1}(U_n)) \subset \phi U \subset U.
\]
This proves the claim.

Now consider the the sequence of maps
\[
	\cdots
	\overset{\alpha}{\longrightarrow} V/U_2
	\overset{\phi}{\longrightarrow}   V/\alpha^{-1}(U_1)  
	\overset{\alpha}{\longrightarrow} V/U_1
	\overset{\phi}{\longrightarrow}   V/\alpha^{-1}(U_0)  
	\overset{\alpha}{\longrightarrow} V/U_0
\]
Since these are injective maps between finite-dimensional vector spaces, all but finitely many
of them are isomorphisms. So for $ n $ large enough the top and bottom
maps in the commutative square
\[
\begin{CD}
	V/\alpha^{-1}(U_n) @>\alpha>> V/U_n \\
	@V{\phi\alpha}VV  @V{\alpha\phi}VV \\
	V/\alpha^{-1}(U_n) @>\alpha>> V/U_n.
\end{CD}
\]
are isomorphisms and it follows that $ \phi\alpha $ and $ \alpha\phi $
have the same characteristic power series.
\end{proof}

\begin{corollary}\label{onepluscommutators}
Let $ N>0 $ be an integer. Assume that all compositions
\[
	\phi, \phi\alpha, \alpha\phi, \phi^2, \cdots
\]
of endomorphisms in $ \{ \phi, \alpha \} $ containing at least one endomorphism
$ \phi $ and at most $ N-1 $ endomorphisms $ \alpha $, are locally contracting. 
Let $ \gamma $ denote $ 1 + \phi $. Then
\[
	\frac{ 1 - \gamma\alpha Z }{ 1 - \alpha\gamma Z } \,\,\mod Z^N
\]
is a nuclear endomorphism of $ V[[Z]]/Z^N $ and
\[
	\det_{k[[Z]]/Z^N}
	\!\left(
		\frac{ 1 - \gamma\alpha Z }{ 1 - \alpha\gamma Z }
		\,\big|\, V
	\right) = 1 .
\]
\end{corollary}

\begin{proof}
We can inductively find uniquely determined endomorphisms $ \phi_n $  so that
\[
	\frac{ 1 - \gamma\alpha Z }
		{ 1 - \alpha Z }
	 = \prod_{n=1}^{\infty} ( 1 - \phi_n \alpha Z^n ) 
\]
and
\[
	\frac{ 1 -\alpha\gamma Z }
		{ 1 - \alpha Z }
	 = \prod_{n=1}^{\infty} ( 1 - \alpha\phi_n Z^n ).
\]
Each $ \phi_n $ is a (non-commutative) polynomial in $ \phi $ and $\alpha $ whose
constituting monomials contain at least one endomorphism
$ \phi $ and at most $ n-1 $ endomorphisms $ \alpha $.
So using the preceding Theorem \ref{commutators} we conclude
\[
	\det_{k[[Z]]/Z^N}
	\!\left(
		\frac{ 1 - \gamma\alpha Z }{ 1 - \alpha\gamma Z }
		\,\big|\, V
	\right) 
	= \prod_{n=1}^{N-1}
	\det_{k[[Z]]/Z^N}
	\!\left(
		\frac{ 1 - \phi_n \alpha Z^n }{ 1 - \alpha \phi_n Z^n }
		\,\big|\, V
	\right) = 1. \qedhere
\]
\end{proof}

\section{A trace formula}\label{sectraceformula}

In this section we establish a trace formula which will enable us to express
$ L( E/R ) $ as the determinant of $ 1 + \Theta $ for some nuclear endomorphism
$ \Theta $ of some compact $ k[[T^{-1}]] $-module.

Similar techniques have been used by
Taguchi and Wan \cite{Taguchi97}, Anderson \cite{Anderson00}, B\"ockle and Pink \cite{Boeckle09},
and by Lafforgue \cite{Lafforgue09} 
to study $L$-functions of $\phi$-sheaves. Our approach differs in that
by allowing more general operators, we are able to compute the $ L $-\emph{value} directly without first computing the $ L $-\emph{function}. Also, we will work directly with the Drinfeld module $ E $, and not with its associated $ \phi $-sheaf or $ t $-motive.

\bigskip

Let $ X $ be a connected projective scheme, smooth of dimension $ 1 $ over $ k $.  
Let $ Y \subset X $ be a non-empty affine open sub-scheme, and $ R $ the ring of
regular functions on $ Y $. Let $ K $ be the function field of $ X $ and denote by
$ K_\infty $ the product of the completions $ K_x $ for all $ x\in X\backslash Y $.
Let $\| \cdot \|\colon K_\infty \to \R$ be the absolute value defined as the maximum of
the normalized absolute values on the $K_x$.

Let $ M $ be a finitely generated projective $ R $-module and
 $ \tau_M\colon M \to M $ a $k$-linear endomorphism satisfying
\[
	\tau_M( rm ) = r^q \tau_M( m )
\]
for all $ r \in R $ and $ m \in M $. 

Let $ \|\cdot \| $ be a norm for the free $ K_\infty $-module $ M \otimes_R K_\infty $.
This norm induces a norm $ \|\cdot\| $ on the compact $ k $-vector space
$ M \otimes_R \frac{K_\infty}{R} $.

Denote by $ R\{ \tau \} $ the twisted polynomial ring whose elements are polynomials
\[
	r_0 + r_1\tau + \cdots + r_d \tau^d 
\]
in $ \tau $, and where multiplication is defined by the commutation rule
 $ \tau r = r^q\tau $ for all
$ r \in R $. By sending $ \tau $ to $ \tau_M $ we have $ k $-algebra homomorphisms 
\[
	R\{\tau\} \to \End_k( \star ) 
\]
 with $ \star $ either $ M $, $ M\otimes_R \frac{K_\infty}{R} $, or $ M/\m M $ for some
 maximal ideal $ \m \subset R $.

\begin{proposition}\label{frobcontract}
For every $ \phi \in R\{\tau\}\tau $ the induced endomorphism
\[
	\phi\colon M \otimes_R \frac{K_\infty}{R} \to M \otimes_R \frac{K_\infty}{R}
\]
is locally contracting.
\end{proposition}

\begin{proof} Clear. In fact any $ c $  with $ 0 < c < 1 $ will do.
\end{proof}

Denote by $ R\{ \tau \}[[ Z ]] $ the ring of formal power series 
\[
	\phi_0 + \phi_1 Z + \phi_2 Z^2 + \cdots 
\]
with $\phi_i \in R\{\tau\}$, where the variable $ Z $ is central. By 
 Proposition \ref{frobcontract} any $ \Theta \in R\{ \tau \}[[Z]]\tau $
defines a nuclear endomorphism of $ \left( M \otimes_R \frac{K_\infty}{R} \right)\![[Z]] $.

The main goal of this section is to establish the following theorem.

\begin{theorem}[Trace formula]\label{traceformula}
For all $ \Theta \in R\{ \tau \}[[Z]] \tau Z $ the infinite product
\[
	 \prod_{\m} \det_{k[[Z]]}\!\left(\, 1 + \Theta \,\big|\, M/\m M \,\right)
\]
converges to
\[
	\det_{k[[Z]]}\!\left( 1 + \Theta \,\big|\,
		 M \otimes_R \frac{K_\infty}{R}\right)^{-1}
\]
in $ k[[Z]] $.
\end{theorem}

This theorem is a variation on Anderson's trace formula, see
Remark \ref{andersontraceformula} at the end of this section. 

We will first show that the validity of Theorem \ref{traceformula}  is invariant under localization. Let $ \p $ be a maximal ideal in $ R $ and denote by
$ R[\p^{-1}] $ and $ M[\p^{-1}] $ the localizations of $ R $ and $ M $ at the multiplicative
system $ \p\cup\{1\} \subset R $.  The endomorphism  
$ \Theta $ of $M[[Z]]$ extends naturally to an endomorphism of
$ M[\p^{-1}][[Z]] $. Let $ K_\p $ be the completion of $ K $ at $ \p $.

\begin{lemma}[Localization]\label{localization}
For all $ \Theta $ as in Theorem \ref{traceformula} we have
\begin{equation}\label{local1}
	 \det_{k[[Z]]}\!\left( 1 + \Theta \,\big|\, M/\p M \right) =
	 \frac{ \det_{k[[Z]]}\!\left( 1 + \Theta \,\big|\,
		M[\p^{-1}] \otimes_{R[\p^{-1}]} \frac{K_\infty\times K_\p}{R[\p^{-1}]}\right) }
		{ \det_{k[[Z]]}\!\left( 1 + \Theta \,\big|\,
		M \otimes_R \frac{K_\infty}{R}\right) }.
\end{equation}
\end{lemma}

\begin{proof}
Denote the completion of $ R $ at $ \p $ by $ R_{\p} $. Since $ K_\p = R_\p + R[\p^{-1}] $
and $ R = R_\p \cap R[\p^{-1}] $ we have a natural short exact sequence
\[
	0 \to R_\p \to \frac{ K_\infty \times K_\p }{ R[\p^{-1}] }
	\to \frac{ K_\infty }{ R } \to 0
\]
inducing a short exact sequence
\[
	0 \to M \otimes_R R_\p \to
	M[\p^{-1}] \otimes_{R[\p^{-1}]} \frac{ K_\infty \times K_\p }{ R[\p^{-1}] }
	\to M\otimes_R \frac{ K_\infty }{ R } \to 0
\]
which is respected by all coefficients $ \theta_n $ of $ \Theta $. It follows that
\[
	 \det_{k[[Z]]}\!\left( 1 + \Theta \,\big|\, M \otimes_R R_\p \right) =
	 \frac{ \det_{k[[Z]]}\!\left( 1 + \Theta \,\big|\,
		M[\p^{-1}] \otimes_{R[\p^{-1}]} \frac{K_\infty\times K_\p}{R[\p^{-1}]}\right) }
		{ \det_{k[[Z]]}\!\left( 1 + \Theta \,\big|\,
		M \otimes_R \frac{K_\infty}{R}\right) }.
\]
But since the subspace $ M \otimes_R \p R_\p \subset M\otimes_R R_\p $ is a common nucleus
for the coefficients $ \theta_n $ of $ \Theta $ we also have
\[
	\det_{k[[Z]]}\!\left( 1 + \Theta \,\big|\, M \otimes_R R_\p \right) =
	\det_{k[[Z]]}\!\left( 1 + \Theta \,\big|\, M/\p M \right),
\]
which proves the lemma.
\end{proof}

\begin{proof}[Proof of Theorem \ref{traceformula}]
Fix positive integers $ D $ and $ N $ and consider the subset
$ \mathcal{S}_{D,N} \subset R\{\tau\}[[Z]]/Z^N $ defined by
\[
	\mathcal{S}_{D,N} := 
	\left\{ 1 + \sum_{n=1}^{N-1} \phi_n Z^n\, \big|\,
		\deg( \phi_n ) < \frac{nD}{N} \text{ for all } n< N \right\}.
\]
The set $ \mathcal{S}_{D,N} $ is a group under multiplication.

We claim that for all $ 1 + \Theta \in \mathcal{S}_{D,N} $ the product
\[
	 \prod_{\m} \det_{k[[Z]]/Z^N}\!\left( 1 + \Theta \,\big|\, M/\m M \right)
\]
converges to
\[
	\det_{k[[Z]]/Z^N}\!\left( 1 + \Theta \,\big|\,
		 M \otimes_R \frac{K_\infty}{R}\right)^{-1}.
\]
Since $ D $ and $ N $ are arbitrary, this claim implies Theorem \ref{traceformula}.

We now apply a variation on a trick of Anderson \cite[Prop. 9]{Anderson00}. 
By Lemma \ref{localization} we may and will assume
that $ R $ has no residue fields of degree $ d < D $ over $ k $. For every $ d<D $ 
we can then find a finite collection of
$ f_{dj} $ and $ a_{dj} $ in $ R $ so that
\[
	1 = \sum_{j} f_{dj} ( a_{dj}^{q^d} - a_{dj} ).
\]
In particular, for every $ r \in R $, $ n<N $ and $ d<D $ we have
\[
	1 - r \tau^d Z^n \equiv
	\prod_j \frac{ 1 - (rf_{dj}\tau^d)a_{dj}Z^n }{ 1 - a_{dj}(rf_{dj}\tau^d)Z^n }
	\,\,\mod Z^{n+1}.
\]
Using this congruence, it is easy to show that the group $ \mathcal{S}_{N,D} $
is generated by elements of the form
\[
	\frac{ 1 - (s\tau^d)aZ^n }{ 1 - a(s\tau^d)Z^n }
\]
with $ a, s \in R $ and $ d $ and $ n $ positive integers. For these elements 
the claim holds thanks to Theorem \ref{commutators},
and by the  multiplicativity of determinants (Proposition \ref{endmultiplicativity})
we may conclude that it holds for all $ 1+\Theta $ in
 $ \mathcal{S}_{D,N} $.
\end{proof}

\begin{remark} \label{andersontraceformula}
The formulation and proof of Theorem \ref{traceformula} grew out of an attempt to
understand and simplify the proof of Anderson's trace formula \cite[Thm 1]{Anderson00}. We
sketch how to deduce Anderson's result from Theorem \ref{traceformula}.

Let $ \Omega_{R/k} $ denote the $ R $-module of K\"ahler differentials of $ R $ over $ k $.
We have a pairing
\[
	\frac{ K_\infty }{ R } \times \Omega_{R/k} \to k
\]
given by
\[
	( f, \omega ) \mapsto \sum_{x\in X\backslash Y} \Tr_{k(x)/k} \Res_x ( f\omega ).
\]
This pairing is perfect in the sense that it identifies either factor with the space
of continuous linear forms on the other. The transpose of the continuous endomorphism 
\[
	\tau\colon \frac{ K_\infty }{ R } \to \frac{ K_\infty }{ R }\colon
	x \to x^q
\]
is the $q$-Cartier operator $ c\colon \Omega_{R/k} \to \Omega_{R/k} $.
Similarly, for any finitely
generated projective $ R $-module $ M $ and any $ \tau_M\colon M \to M $ as above the
induced homomorphism 
\[
	\tau_M\colon M\otimes_R \frac{ K_\infty }{ R }  
	\to M \otimes_R \frac{ K_\infty }{ R } 
\]
has a ``Cartier-linear'' transpose endomorphism
\[
	 c_M\colon \Hom_R( M, \Omega_{R/k} )  \to \Hom_R( M, \Omega_{R/k} ).
\]
Using this one deduces Anderson's trace formula from Theorem \ref{traceformula} applied
to $ \Theta = 1 - \tau_M Z $. 
\end{remark}

\section{Ratio of co-volumes} \label{secratio}

Let $ V $ be a finite-dimensional $ k((t^{-1})) $-vector space. For
lattices $ \Lambda_1 $ and $ \Lambda_2 $ in $ V $ we have defined in 
\S \ref{intro}, Proposition \ref{defreg} 
the quantity $  [\Lambda_1:\Lambda_2] \in  k((T^{-1})) $.
In the present section we will express $ [\Lambda_1 : \Lambda_2] $ as a determinant
in the sense of \S  \ref{secdeterminants}.

\bigskip

Let $ \|\cdot \|\colon V \to \R_{\geq 0} $ be a norm on the
$k((t^{-1}))$-vector space $ V $ and let $ \Lambda_1 $ and $ \Lambda_2 $ be lattices in $ V $.

\begin{definition}\label{deftangent}\label{definftangent}
Let $N$ be a positive integer.
We say that a continuous $k$-linear map 
\[
	\gamma\colon V/\Lambda_1 \to V/\Lambda_2
\]
is \emph{$N$-tangent to the identity} (on $ V $) if there is an open $ U \subset V $
such that
\begin{enumerate}
	\item $ U \cap \Lambda_1 = U \cap \Lambda_2 = 0 $;
	\item $ \gamma $ restricts to an isometry between the images of $ U $;
	\item for all $ v \in U $ we have
		$ \| \gamma( v ) - v \| \leq \| t^{-N} \| \cdot \| v \| $.
\end{enumerate}
\end{definition}

\begin{definition}
We say that $ \gamma $
is \emph{infinitely tangent to the identity} (on $ V $) if it is $N$-tangent 
to the identity for every positive integer $ N $. 
\end{definition}

\begin{proposition}\label{powerseries}
Let $F$ be a finite extension of $k((t^{-1}))$ and let $\| \cdot \|$ be an 
extension of the absolute value on $\| \cdot \|$ on $k((t^{-1}))$. Let $\Lambda_1$
and $\Lambda_2$ be lattices in the $k((t^{-1}))$-vector space $F$. 
Let $\gamma\colon F\to F$ be a $k$-linear map defined by an everywhere convergent power series
\[
	\gamma(x) = x + \gamma_1 x^q + \gamma_2 x^{q^2} + \cdots,
\]
and such that $\gamma(\Lambda_1)\subset \Lambda_2$. Then the induced map
\[
	\gamma\colon F/\Lambda_1 \to F/\Lambda_2
\]
is infinitely tangent to the identity on $F$.
\end{proposition}

\begin{proof}
Let $N$ be a positive integer.
By the convergence of the power series the coefficients $\gamma_i$ are bounded and so
there exists an $\epsilon>0$ such that
\[
	\| \gamma(x) - x \|  < \|t^{-N}\| \cdot \| x \|
\]
for all $x\in F$ with $\|x\| < \epsilon$. Shrinking $\epsilon$ if necessary, we may assume that
$\| \lambda \| \geq \epsilon$ for all nonzero $\lambda$ in $\Lambda_1$ or $\Lambda_2$. Then
by the non-archimedean inverse function theorem \cite[2.2]{Igusa00} the open subset
$U = \{ x \in F \colon \| x \| < \epsilon \} $ satisfies the requirements. 
\end{proof}

Let $ H_1 $ and $ H_2 $ be finite $ k[t] $-modules and denote $ V/\Lambda_i \times H_i $
by $ M_i $. Let
\[
	 \gamma\colon M_1 \to M_2
\]
be a continuous $k$-linear map. We say that $ \gamma $ is $N$-tangent to the identity 
on $V$ if the composition
\[
	V/\Lambda_1 \into M_1 \overset{\gamma}{\longrightarrow} M_2 \onto V/\Lambda_2
\]
is $N$-tangent to the identity, and we say that $\gamma$ is infinitely tangent to the identity if it is $N$-tangent to the identity on $V$ for all positive integers $N$.

For every $k$-linear isomorphism $ \gamma\colon M_1 \to M_2 $ we define an endomorphism 
endomorphism 
\[
	\Delta_\gamma = \frac{1-(\gamma^{-1}t\gamma)T^{-1}}{1- tT^{-1}} - 1   
\]
of $ M_1[[T^{-1}]] $.  We have $ \Delta_\gamma = \sum_{n=1}^\infty \delta_n T^{-n} $ with
\[
	\delta_n = (t - \gamma^{-1}t\gamma)t^{n-1}.
\]

\begin{theorem} \label{detreg}
If $ \gamma $ is infinitely tangent to the identity on $V$ then
$ \Delta_\gamma $ is nuclear and
\begin{equation}\label{eqdetreg}
	\det_{k[[T^{-1}]]}\!\left( 1 + \Delta_\gamma \,\big|\, M_1 \right)
	= [ \Lambda_1 : \Lambda_2 ] \frac{|H_2|}{|H_1|}.	
\end{equation}
\end{theorem}

\begin{remark}\label{volumes}
It is natural to think of 
\[
	 [ \Lambda_1 : \Lambda_2 ] \frac{|H_2|}{|H_1|}
\]
as the ratio of the ``volumes'' of
$ M_2 $ and $ M_1 $. It seems plausible that one could work out a theory
of $ k((T^{-1})) $-valued ``Haar measures'' on topological $ k[t] $-modules in which 
such a statement would have a precise meaning. (See also Remark \ref{tamagawa}.)
\end{remark}

\begin{remark}
If $ V=0 $ then $ M_1 =H_1 $ and $ M_2=H_2 $ are finite and the theorem 
trivially holds since
\[
	\det_{k[[T^{-1}]]}\!\left( 1 + \Delta_\gamma \,\big|\, M_1 \right)
	= \frac{ \det_{k[[T^{-1}]]}\!\left( 1 - tT^{-1} \,\big|\, H_2 \right) }
		{ \det_{k[[T^{-1}]]}\!\left( 1 - tT^{-1} \,\big|\, H_1 \right) }
	= \frac{|H_2|}{|H_1|}.
\]
\end{remark}

\begin{remark}
By Theorem \ref{detreg} the existence of a $k$-linear
isomorphism $ \gamma\colon M_1 \to M_2 $ which is infinitely tangent to the identity 
implies that $ [ \Lambda_1 : \Lambda_2 ] \frac{|H_2|}{|H_1|} $,
which is \emph{a priori} an element of $ k((T^{-1}))^\times $,
has valuation zero. It is not hard to show that the converse is also true.
\end{remark}

The proof of Theorem \ref{detreg} makes up the rest of this section. 

\begin{lemma} 
If $ \gamma $ is $N$-tangent to the identity then
$ \Delta_\gamma $ modulo $ T^{-N} $ is a nuclear endomorphism of
$ M_1[[T^{-1}]]/T^{-N} $.
\end{lemma}

\begin{proof}
For all $ n<N $ and for all $ m\in M_1 $ sufficiently small we have
\begin{eqnarray*}
	\| \delta_n(m) \|  &=& 
	\| \gamma^{-1}\big( \gamma(t^nm)-t^nm \big)
	+ \gamma^{-1}\big( t^nm-t\gamma(t^{n-1}m) \big) \| \\
	&\leq & \| t^{-1} \| \cdot \| m \|,
\end{eqnarray*}
so $ \delta_n $ is locally contracting.
\end{proof}

\begin{lemma}[Independence of $\gamma$]\label{independence}
Assume that $ \gamma_1\colon M_1 \to M_2 $ and $ \gamma_2\colon M_1 \to M_2 $ are
$N$-tangent to the identity. Then
\[
	\det_{k[[T^{-1}]]/T^{-N}}\!\left( 1 + \Delta_{\gamma_1} \,\big|\, M_1 \right) =
	\det_{k[[T^{-1}]]/T^{-N}}\!\left( 1 + \Delta_{\gamma_2} \,\big|\, M_1 \right)
\]
in $ k[[T^{-1}]]/T^{-N} $.
\end{lemma}

\begin{proof}
By Proposition \ref{endmultiplicativity} it suffices to show  that
\[
	\det_{k[[T^{-1}]]/T^{-N}}\!\left( 1 + \Delta_{\gamma} \,\big|\, M \right) = 1
\]
when $ \gamma $ is an automorphism of $ M = M_1 $ which is $N$-tangent to the identity.
This follows from Corollary \ref{onepluscommutators}, with $\alpha=t\gamma^{-1}$
and $ \phi = \gamma - 1 $.
\end{proof}

\begin{lemma}[Torsion]\label{torsion}
If both $ \Lambda_1 $ and $ \Lambda_2 $ are contained in a common
$ k[t] $-lattice $ \Lambda$, then Theorem \ref{detreg} holds.
\end{lemma}

\begin{proof}
By the previous lemma we are free to choose $ \gamma $.
Let $ U \subset V $ be an open subspace that is a complement of $ \Lambda \subset V $.
We have decompositions
\[
	M_i = U \times \Lambda/\Lambda_i \times H_i \quad\quad (i=1,2.)
\]
Choose a $ k $-isomorphism 
\[
	\gamma: M_1  \to M_2 
\]
that is the identity on $ U $ and that maps $ \Lambda/\Lambda_1 \times H_1 $  isomorphically
to $ \Lambda/\Lambda_2 \times H_2 $. Clearly $ \gamma $ is tangent to the identity on $ V $.
The finite subspace 
\[
	\Lambda/\Lambda_1\times H_1 \subset M_1
\]
is preserved under the coefficients $ \delta_n $ of $ \Delta_\gamma $ and they induce the
zero map on the quotient, hence 
\begin{eqnarray*}
	\det_{k[[T^{-1}]]}\!\left( 1 + \Delta_\gamma \,\big|\, M_1 \right)
	&=&
	\det_{k[[T^{-1}]]}\!\left(
		1 + \Delta_{\gamma}\,\big|\, \Lambda/\Lambda_1 \times H_1 \right) \\
	&=&
	\frac{|\Lambda/\Lambda_2 \times H_2|}{|\Lambda/\Lambda_1 \times H_1|} \\
	&=& [ \Lambda_1 : \Lambda_2 ] \frac{|H_2|}{|H_1|},
\end{eqnarray*}
as desired.
\end{proof}

\begin{lemma}[Approximation]\label{approx}
Let $ N $ be a positive integer, $ U \subset V $ an open sub-$k[[t^{-1}]]$-module and
$ \sigma $ an element of $ \GL_{k[[t^{-1}]]}(U) $ inducing the identity on $ U / t^{-N} U $.
Then for every lattice $ \Lambda $ in $ V $ the isomorphism
\[
	\sigma \colon V/\Lambda \to V/\sigma(\Lambda)
\]
is $N$-tangent to the identity and we have
\[
	\det_{k[[T^{-1}]]/T^{-N}}\!\left( 1 + \Delta_\sigma \,\big|\, V/\Lambda \right)
	=
	1 \mod{ T^{-N} }.
\]
\end{lemma}

\begin{proof}
Clearly $ \sigma $ is $N$-tangent to the identity, and since $ \sigma $ is $ k[t] $-linear
we have $ \Delta_\sigma = 0 $.
\end{proof}

\begin{proof}[Proof of Theorem \ref{detreg}]
Let $ N $ be an arbitrary positive integer. It suffices to prove (\ref{eqdetreg}) modulo
$ T^{-N} $. But by Lemma \ref{approx} and Lemma \ref{independence}
we may replace $ \Lambda_1 $ with a lattice $ \Lambda_1' $
which is contained in a common over-lattice with $ \Lambda_2 $, and such that 
\[
	\frac{ [ \Lambda_1 : \Lambda_2 ] }{ [ \Lambda_1' : \Lambda_2 ] }
	\in 1 + T^{-N}k[[T^{-1}]]
\]
and 
\[
	\frac{
		\det_{k[[T^{-1}]]}\!\left( 1 + \Delta_\gamma \,\big|\, M_1 \right) }{
		\det_{k[[T^{-1}]]}\!\left( 1 + \Delta_{\gamma'} \,\big|\, M_1' \right)
	} \in 1 + T^{-N}k[[T^{-1}]]
\]
for every $ \gamma'\colon M_1' \to M_1 $ infinitely tangent to the identity. Together
with Lemma \ref{torsion} this proves the theorem modulo $ T^{-N} $.
\end{proof}

\section{Proof of the main result}

We now have at our disposal all the ingredients 
necessary to prove the main result, Theorem \ref{mainthm}. 

We recall the set-up. The ring $ R $ is
the integral closure of $ k[t] $ in a finite extension $ K $ of $ k(t) $, 
and $ K_\infty $ denotes the $ k((t^{-1})) $-algebra $ R \otimes_{k[t]} k((t^{-1})) $.
Let $ E $ be a Drinfeld module over $ R $ defined by
\[
	\phi_E(t) = t + r_1 \tau + \cdots + r_n \tau^n \in R\{\tau\}
\]
where $ \tau $ is the $ q $-th power Frobenius endomorphism of the additive group.

Consider the power series
\[
	\Theta := \frac{ 1-\phi_E(t)T^{-1} }{ 1 - tT^{-1} } - 1 = 
	\sum_{n=1}^\infty (t-\phi_E(t))t^{n-1}T^{-n} 
\]
in $ R\{\tau\}[[T^{-1}]]\tau T^{-1} $.
We can express the special value $ L( E/R ) $ using $ \Theta $ as follows:
\begin{eqnarray*}
	L( E/R ) &=& \prod_{\m} \frac{ |R/\m | }{ |E(R/\m)| } \\
		&=&
	 \prod_\m \det_{k[[T^{-1}]]}\!\left( 1 + \Theta \,\big|\, R/\m \right)^{-1}.
\end{eqnarray*}
Theorem \ref{traceformula} applied to $\Theta$ acting on $ R[[T^{-1}]] $ then gives
\[ 
	L( E/R ) =
	\det_{k[[T^{-1}]]}\!\left( 1 + \Theta \,\big|\, \frac{ K_\infty }{ R } \right).
\]
Now consider the short exact sequence of $ k[t] $-modules
\[
	0 \longrightarrow 
	\frac{ K_\infty }{ \exp_E^{-1} E(R)  }
	\overset{\exp_E}{\longrightarrow}
	\frac{E(K_\infty)}{E(R)} \longrightarrow
	H(E/R) \longrightarrow 0.
\]
Since the  $ k[t] $-module on the left is divisible this
sequence splits. The choice of a section gives an isomorphism
\[
	 \frac{K_\infty}{\exp^{-1} E(R)} \times H(E/R) \overset\sim\to
	 \frac{E(K_\infty)}{E(R)}.
\]
Denote by $ \gamma $ the composition 
\[
	\frac{K_\infty}{\exp_E^{-1}(E(R))} \times H(E/R) \overset\sim\to
	 \frac{E(K_\infty)}{E(R)} \overset\sim\to \frac{K_\infty}{R}
\]
with the tautological map. Then $\gamma$ is a $k$-linear isomorphism (but in general it is not $k[t]$-linear).

Let $\|\cdot\|\colon K_\infty \to \R$ be an absolute value extending the given absolute value
on $k((t^{-1}))$.

\begin{lemma}\label{exptangent}
$ \gamma $ is infinitely tangent to the identity on $ K_\infty $.
\end{lemma}

\begin{proof} This boils down to the statement that
\[
	\exp_E\colon \frac{K_\infty}{\exp^{-1} E(R)} \to \frac{K_\infty}{R}
\]
is infinitely tangent to the identity, which follows from Proposition 
\ref{powerseries}.
\end{proof}

Since we have an equality
\[
	1 + \Theta = \frac{ 1 - \gamma t\gamma^{-1}T^{-1} }{ 1 - tT^{-1} }  
\]
of endomorphisms of $ \frac{K_\infty}{R}[[T^{-1}]] $ we  conclude using Theorem \ref{detreg}:
\begin{eqnarray*}
	L( E/R ) &=& 
	\det_{k[[T^{-1}]]}\!\left( 1 + \Theta \,\big|\, \frac{ K_\infty }{ R } \right) \\
	&=&
	\left[\, R : \exp_E^{-1} E(R) \,\right] \cdot \left| H(E/R) \right|.
\end{eqnarray*}

\section*{Acknowledgements}

The author is grateful to David Goss for his numerous comments and his enthusiastic encouragement,
and to the anonymous referees for several useful suggestions. Parts of the research leading to this paper were carried out at the University of Pennsylvania and at the 
\'Ecole Polytechnique F\'ed\'erale de Lausanne. The author is supported by a VENI grant of the Netherlands Organisation for Scientific Research (NWO). 

\small
\bibliographystyle{plain}
\bibliography{../../master}

\end{document}